\documentclass[12pt]{amsart}
\usepackage{amsmath, amssymb, natbib}

\DeclareMathAlphabet{\curly}{U}{rsfs}{m}{n}

\newtheorem{thm}{Theorem}
\newtheorem{cor}{Corollary}

\newcommand{\Vol}{\operatorname{Vol}}   

\newcommand{\g}{\ensuremath{\gamma}}

\newcommand{\lam}{\ensuremath{\lambda}}
\newcommand{\eps}{\ensuremath{\varepsilon}}

\newcommand{\bxi}{\ensuremath{\boldsymbol{\xi}}}

\newcommand{\fl}[1]{{\ensuremath{\left\lfloor {#1} \right\rfloor}}}

\newcommand{\pfrac}[2]{{\left(\frac{#1}{#2}\right)}}

\newcommand{\be}{\begin{equation}}
\newcommand{\ee}{\end{equation}}
\newcommand{\benn}{\begin{equation*}}   
\newcommand{\eenn}{\end{equation*}}

\renewcommand{\b}{\ensuremath{\beta}}
\renewcommand{\a}{\ensuremath{\alpha}}

\renewcommand{\(}{\left(}
\renewcommand{\)}{\right)}

\begin{document}

\title{Generalized Smirnov statistics and the distribution of prime factors}
\author{Kevin Ford}

\dedicatory{Dedicated to Jean-Marc Deshouillers on the occasion of his 60th
  birthday} 

\address{Department of Mathematics, 1409 West Green Street, University
of Illinois at Urbana-Champaign, Urbana, IL 61801, USA}
\email{ford@math.uiuc.edu}

\date{\today}
\thanks{Research supported by National Science Foundation grant
DMS-0555367.}

\begin{abstract} We apply recent bounds of the author for generalized Smirnov
statistics to the distribution of integers whose prime factors satisfy
certain systems of inequalities. 
\end{abstract}

\maketitle

%
%
%
\section{Introduction}\label{sec:intro}
%
%
%

For a positive integer $n$, denote by $p_1<p_2<\cdots<p_{\omega(n)}$ the
sequence of distinct prime factors of $n$.  In this note, we study integers
for which
\be\label{pjlower}
\log_2 p_j \ge \a j - \b \quad (1\le j\le \omega(n))
\ee
or
\be\label{pjupper}
\log_2 p_j \le \a j + \b \quad (1\le j\le \omega(n)),
\ee
where  $\a\ge 0$ and $\log_2 y$ denotes $\log\log y$.  
The distribution of integers satisfying \eqref{pjlower} is important in
the study of the distribution of divisors of integers (see \cite{F}; Ch. 2 of
\cite{Divisors}).  We present here estimates for
\begin{align*}
N_k(x;\a,\b) &= \# \{ n\le x : \omega(n)=k, \eqref{pjlower} \}, \\
M_k(x;\a,\b) &= \# \{ n\le x : \omega(n)=k, \eqref{pjupper} \}.
\end{align*}

It is a relatively simple matter, at least heuristically, to reduce the
estimation of $N_k(x;\a,\b)$ and $M_k(x;\a,\b)$ to the estimation of
a certain probability connected to Kolmogorov-Smirnov statistics.  
Let us focus on the upper bound for $N_k(x;\a,\b)$.  If we suppose
that $p_k \ge x^c$ for some small $c$, then for each choice of
$(p_1,\ldots,p_{k-1})$, the number of possible $p_k$ is 
$\ll x/(p_1\cdots p_{k-1} \log x)$.  
Since $\sum_{p\le y} 1/p \approx \log_2 y$,
given a well-behaved function $f$, by partial summation we anticipate that
\be\label{approx}
\sum_{p_1<\cdots <p_{k-1} \le x} \frac{f\( \frac{\log_2 p_1}{\log_2 x},
  \cdots, \frac{\log_2 p_{k-1}}{\log_2 x} \)}{p_1\cdots p_{k-1}} \approx
(\log_2 x)^{k-1} \idotsint\limits_{0\le \xi_1 \le \cdots\le \xi_{k-1}\le 1} 
f(\bxi)\, d\bxi,
\ee
where $\bxi=(\xi_1,\ldots,\xi_{k-1})$.

Let $U_1, \ldots, U_{m}$ be independent, uniformly distributed random
variables in $[0,1]$ and let $\xi_1,\ldots,\xi_{m}$ be their 
order statistics ($\xi_1$ is the smallest of the $U_i$, $\xi_2$ is the
next smallest, etc.).  Taking $m=k-1$, 
the right side of \eqref{approx} is equal
to $(\log_2 x)^{k-1}/(k-1)!$ times the expectation of
$f(\xi_1,\ldots,\xi_{k-1})$.  Letting $f$ be 1 if \eqref{pjlower} holds and 0
otherwise,  the expectation of $f$ is
the probability that $\xi_j \ge (\a j -\b)/\log_2 x$ for each $j$.

In general, let $Q_m(u,v)$ be
the probability that $\xi_i \ge \frac{i-u}{v}$ for $1\le i \le m$.
Equivalently, if $u\ge 0$ then
$$
Q_m(u,v) = \text{Prob} \( F_m(t) \le \frac{vt+u}{m} \;\; (0\le t\le 1) \), 
$$
where $F_m(t)=\frac{1}{m}\sum_{U_i\le t} 1$ is the associated
 empirical distribution function. 
The first estimates for $Q_m(u,v)$ were given in 1939 by N. V. Smirnov
\cite{Sm1}, who 
proved for each fixed $\lam\ge 0$ the asymptotic formula
\be\label{Smirnov}
Q_m(\lam \sqrt{m},m) \to 1 - e^{-2\lam^2} 
\qquad (m\to\infty).
\ee
The sharpest and most general bounds are due to the author \cite{smir}; see
also \cite{kol2}.
For convenience, write $w=u+v-m$.  
Uniformly in $u>0$, $w>0$ and $m\ge 1$, we have
\be\label{Qmuv}
Q_m(u,v) = 1 - e^{-2 uw/m} + O\pfrac{u + w}{m}.
\ee
Moreover, 
\be\label{Q}
Q_m(u,v) \asymp \min \( 1 , \frac{uw}{m}\) \qquad (u\ge 1, w\ge 1).
\ee
See \cite{smir} for more information about the history of such bounds and
techniques for proving them.  A short proof of weaker bounds is given in
\S 11 of \cite{F}.

Returning to our heuristic estimation of $N_k(x)$ (and assuming that
a similar lower bound holds), we find that
$$
N_k(x) \approx \frac{x (\log_2 x)^{k-1}}{(k-1)! \log x} 
Q_{k-1} \( \frac{\b}{\a}, \frac{\log_2 x}{\a} \).
$$
We have (cf. Theorem 4 in \S II.6.1 of \cite{Tenbook})
\be\label{pik}
\pi_k(x) := \# \{ n\le x : \omega(n)=k \} \asymp_A 
\frac{x (\log_2 x)^{k-1}}{(k-1)! \log x} 
\ee
uniformly for $1\le k\le A\log_2 x$, $A$ being any fixed positive constant.
Thus, we anticipate that
$$
N_k(x;\a,\b) \asymp Q_{k-1} \( \frac{\b}{\a}, \frac{\log_2 x}{\a} \)  \pi_k(x).
$$
Observing that the vectors $(\xi_1,\ldots,\xi_m)$ and
$(1-\xi_m, 1-\xi_{m-1},\ldots,1-\xi_1)$ have identical distributions,
we have
$$
Q_m(u,v) = \text{Prob} \( \xi_i \le \frac{u+v-m-1+i}{v} \;\; (1\le i\le m) \).
$$
Hence, we likewise anticipate that
$$
M_k(x;\a,\b) 
\asymp Q_{k-1} \(  k+\frac{\b-\log_2 x}{\a}, \frac{\log_2 x}{\a} \)
\pi_k(x).
$$

To make our heuristics rigorous, we must impose some 
conditions on $\a$ and $\b$ to ensure among other things that there
are integers satisfying \eqref{pjlower} or \eqref{pjupper}.
To that end, we set
\be\label{Nkuvw}
u= \frac{\b}{\a}, \qquad v = \frac{\log_2 x}{\a}, \qquad
w = u+v - (k-1) = \frac{\log_2 x+\b}{\a} - k +1
\ee
for the estimation of $N_k(x;\a,\b)$ and
\be\label{Mkuvw}
u = k + \frac{\b-\log_2 x}{\a}, \qquad v = \frac{\log_2 x}{\a}, \qquad
w = u + v - (k-1) = \frac{\b}{\a} + 1
\ee
for the estimation of $M_k(x;\a,\b)$.

\begin{thm}\label{Nk}  Suppose $\eps>0$, $A\ge 1$ and $1\le k\le A\log_2 x$.
Assume \eqref{Nkuvw}, $\b\ge 0$, $\a-\b \le A$, $w\ge 1+\eps$ and
\be\label{exp1}
e^{\a (w-1)} - e^{\a(w-2)} \ge 1 + \eps.
\ee
Then, for sufficiently large $x$, depending on $\eps$ and $A$,
$$
N_k(x;\a,\b) \asymp_{\eps,A} \min \( 1 , \frac{(u+1)w}{k}\) \pi_k(x),
$$
the implied constants depending only on $\eps$ and $A$.
\end{thm}

\begin{thm}\label{Mk}  Suppose $A\ge 1$ and $1\le k\le A\log_2 x$.
Assume \eqref{Mkuvw}, $u\ge 1$, $w\ge 0$ and that for $1\le j\le k$, there
are at least $j$ primes $\le \exp \exp (\a j+\b)$.
Then, for sufficiently large $x$, depending on $A$,
$$
M_k(x;\a,\b) \asymp_{A} \min \( 1 , \frac{u(w+1)}{k}\) \pi_k(x),
$$
the implied constants depending only on $A$.
\end{thm}

{\bf Remarks.}  Inequality \eqref{exp1} is necessary, since for large
$k$, \eqref{pjlower} implies 
$$
\log n \ge \sum_{j=1}^k \log p_j \ge \sum_{j=1}^k e^{\a j - \b} \approx
\frac{e^{\a k-\b}}{1-e^{-\a}} = \frac{\log x}{e^{\a(w-1)}-e^{\a(w-2)}}.
$$
The condition $\a-\b\le A$ in Theorem \ref{Nk} 
means that there is no significant
restriction on $p_1$.

It is a simple matter to apply 
the estimates for $N_k(x;\a,\b)$ and $M_k(x;\a,\b)$
to problems of the distribution of prime factors of integers where $\omega(n)$
is not fixed.  In the following, let $\omega(n,t)$ be the number of distinct
prime factors of $n$ which are $\le t$.  It is well-known 
(cf. Ch. 1 of \cite{Divisors}) that $\omega(n,t)$ has normal order $\log_2 t$.
We estimate below the likelihood that $\omega(n,t)$ does not stray too far
from $\log_2 t$ in one direction.


\begin{cor}\label{cor}
Uniformly for large $x$ and $0\le \beta \le \sqrt{\log_2 x}$, we have
\be\label{cor1}
\# \{ n\le x : \forall t, 2\le t\le x,\omega(n,t)\le \max(0,\log_2 t + \b) \} 
\asymp \frac{(\beta+1) x}{\sqrt{\log_2 x}}
\ee
and
\be\label{cor2}
\# \{ n\le x : \forall t, 2\le t\le x, \omega(n,t)\ge \log_2 t - \b \} \asymp
\frac{(\beta+1) x}{\sqrt{\log_2 x}}
\ee
\end{cor}

\emph{Proof of Corollary \ref{cor}}.
The quantity of the left side of \eqref{cor1} is $\sum_k N_k(x;1,\b)$.
Here $u=\b$, $v=\log_2 x$ and $w=\log_2 x + \b - k + 1$.  By Theorem \ref{Nk}
and \eqref{pik},
$$
\sum_{\log_2 x - 2\sqrt{\log_2 x} \le k \le \log_2 x - \sqrt{\log_2 x}}
N_k(x;1,\b) \gg \frac{(\b+1) x}{\sqrt{\log_2 x}},
$$
since $\pi_k(x) \asymp x/\sqrt{\log_2 x}$ for $|k-\log_2 x| \le 2\sqrt{\log_2
  x}$.   This proves the lower bound in \eqref{cor1}.  For the upper bound,
we note that if $k > \log_2 x + \b$, then $N_k(x;1,\b)=0$.  Hence, by
Theorem \ref{Nk} and \eqref{pik},
\begin{align*}
\sum_{k} N_k(x;1,\b) &\ll \sum_{k\le \log_2 x+\b-2} \frac{(\b+1)(\log_2 x +
\b-k+1)}{k}\pi_k(x) + \sum_{\log_2 x + \b-2 < k \le \log_2 x + \b} \pi_k(x) \\
&\ll \frac{(\b+1) x}{\sqrt{\log_2 x}}.
\end{align*}
This proves the upper bound in \eqref{cor1}.

The quantity on the left side of \eqref{cor2} is $\sum_k M_k(x;1,\b-1)$.
Here $v=\log_2 x$, $u=\b + k -\log_2 x$ and $w=\b$.  By Theorem \ref{Mk},
$$
\sum_{\log_2 x + \sqrt{\log_2 x} \le k \le \log_2 x + 2\sqrt{\log_2 x}}
M_k(x;1,\b-1) \gg \frac{(\b+1) x}{\sqrt{\log_2 x},}
$$
proving the lower bound in \eqref{cor2}.
Also by Theorem \ref{Mk},
$$
\sum_{\log_2 x - \b + 1 < k \le 10 \log_2 x} M_k(x;1,\b-1) \ll \frac{(\b+1)
  x}{\sqrt{\log_2 x}}.
$$
If $\omega(n)=k > 10\log_2 x$, then the number, $\tau(n)$, of divisors of $n$
satisfies $\tau(n) \ge 2^{\omega(n)} \ge (\log x)^6$.  Since 
$\sum_{n\le x} \tau(n) \sim x\log x$, the number of $n\le x $ with $\omega(n)
> 10 \log_2 x$ is $O(x/\log^5 x)$.
By \eqref{pik}, the number of $n\le x$ with $\log_2 x - \b - 4 < k \le \log_2
x - \b + 1$ is $O(x/\sqrt{\log_2 x})$.  Finally, suppose $k \le \log_2 x - \b
- 4$.  The number of $n\le x$ for which $d^2|n$ for some $d>\log x$ is
$O(x\sum_{d>\log x} 1/d^2) = O(x/\log x)$.  If there is no such $d$, then
by \eqref{pjupper},
$$
\log n \le 2\log_2 x + \sum_{j=1}^k \log p_j \le 2\log_2 x + \sum_{j=1}^k
e^{j+\b-1} \le 2\log_2 x + 2 e^{k+\b-1} \le \tfrac12 \log x,
$$
thus $n\le \sqrt{x}$.  This completes the proof of the upper bound in
\eqref{cor2}.

Our methods for proving Theorems \ref{Nk} and \ref{Mk}
are borrowed from \cite{F}, especially sections 8, 10 and
12 therein.   The tools there are adequate
for making precise the heuristic argument outlined above when the function $f$
is monotonic in each variable, even if $f$ is discontinuous.
We provide details only for Theorem \ref{Nk}.
In lower bound for $M_k(x;\a,\b)$, we may need to fix several of the
smallest prime factors of $n$, but otherwise the details of the proof
of Theorem \ref{Mk} are very similar.

%
\section{Certain partitions of the primes}\label{partitions}
%

We describe in this section certain partitions of the primes which will
be needed in the proof of Theorems \ref{Nk} and \ref{Mk}.  The constructions
are similar to those given in \S 4 and \S 8 of \cite{F}.

Let $\lambda_0=1.9$ and inductively define $\lambda_j$ to be the largest prime
such that
$$
\sum_{\lambda_{j-1} < p \le \lambda_j} \frac{1}{p} \le 1.
$$
In particular, $\lambda_1=3$ and $\lambda_2=109$.  By Mertens' estimate,
$\log_2 \lambda_j = j + O(1)$.  Let $G_j$ be the set of primes in
$(\lambda_{j-1},\lambda_j]$ for $j\ge 1$.  Then there is an absolute 
constant $K$ so that if $p\in G_j$ then $|\log_2 p -j | \le K$.

Next, let $Q \ge e^{10}$ and $\g = 1/\log Q$.  
If $p\le Q$, then $p^\g \le e$, hence
$p^\g \le 1 + (e-1)\g \log p$.
By Merten's estimates,
$$
\sum_{\substack{p\le Q \\ f\ge 1}} \frac{1}{p^{f(1-\g)}} =
O(1) + \sum_{p\le Q} \( \frac{1}{p}  + (e-1)\g \frac{\log p}{p} \) =
\log_2 Q + O(1).
$$
It follows for an absolute constant $K'$, independent of $Q$, that
the set of primes $p\le Q$ may be partitioned into at most 
$\frac12 \log_2 Q + K'$ sets $E_j$ so that (i) for each $j$,
$$
\sum_{\substack{p\in E_j \\ f\ge 1}} \frac{1}{p^{f(1-\g)}} \le 2
$$
and (ii) for $p\in E_j$,  $|\log_2 p - 2j| \le K'$.
We stipulate that the above sum is $\le 2$ rather than $\le 1$ in order
to accomodate the prime 2.

%
\section{Proof of Theorem \ref{Nk} upper bound}
%

Without loss of generality, suppose that $k$ is large, $(u+1)w\le k/10$,
and $n\ge x/\log x$.
We have $v\le 1.1k$ and
consequently $\a \ge 1/(1.1A)$.  Also, by \eqref{pjlower},
$$
\log_2 p_k \ge \a k - \b = \frac{k-u}{v} \log_2 x \ge \frac{9}{11} \log_2 x.
$$
We may suppose $p_k^2 \nmid n$, as the number of $n\le x$ with
$p_k^2|n$ is $O(x \exp(-(\log x)^{9/11})) = O(\pi_k(x)/k)$.
For brevity, write $x_\ell = x^{1/e^{\ell}}$.  For some integer $\ell$
satisfying $\ell\ge 0$ and $\exp\exp (\a k- \b) \le x_\ell$, we have
$x_{\ell+1} < p_k \le x_{\ell}$.  With $\ell$ fixed, given
$p_1,\ldots,p_{k-1}$ with exponents $f_1,\ldots,f_{k-1}$, 
the number of possibilities for $p_k$ is
$$
\ll \frac{x}{p_1^{f_1} \cdots p_{k-1}^{f_{k-1}} \log x_\ell}
\ll \frac{x^{1-\g/2} e^{\ell}}{(p_1^{f_1} 
\cdots p_{k-1}^{f_{k-1}})^{1-\g}\log x},
$$ 
where $\g=1/\log x_\ell$.  This follows
for $\ell \ge 1$ from $p_1^{f_1} \cdots p_{k-1}^{f_{k-1}} \ge x/(p_k \log x) >
x^{1/2}$.  We conclude that
\be\label{Nku1}
N_k(x;\a,\b) \ll \frac{x}{\log x} \sum_{\ell} e^{\ell-\frac12 e^{\ell}} 
\sum_{\substack{p_1 <\cdots <p_{k-1}\le x_\ell \\ f_1,\ldots,f_{k-1} \ge 1 \\
\eqref{pjlower}}} \frac{1}{(p_1^{f_1}  \cdots p_{k-1}^{f_{k-1}})^{1-\g}}.
\ee

Consider the intervals $E_j$ defined in the previous section corresponding to
$Q=x_\ell$.  Put $J=\fl{\frac12\log_2 x_\ell+K'}$ and define
$j_1,\ldots,j_{k-1}$ by $p_i \in E_{j_i}$.  Let $\mathcal{J}$ denote the set of
tuples $(j_1,\ldots,j_{k-1})$ so that $1\le j_1 \le \cdots \le j_{k-1}\le J$
and such that $j_i\ge \frac12 (\a i - \b - K'-A)$ for every $i$.
Given $p_1,\ldots,p_{k-1}$, let $b_j$ be the number of $p_i$ in $E_j$, for
$1\le j\le J$.  The contribution to the inner sum of \eqref{Nku1} from
those tuple of primes with a fixed $(j_1,\ldots,j_{k-1})$ is
\begin{align*}
&\le \prod_{j=1}^J \frac{1}{b_j!} \biggl( \sum_{p\in E_j, f\ge 1} \frac{1}
{p^{f(1-\g)}} \biggr)^{b_j} \\ 
&\le \frac{2^{k-1}}{b_1! \cdots b_J!}.
\end{align*}
We observe that $1/(b_1!\cdots b_J!)$ is the volume of the region
$(y_1,\cdots,y_{k-1}) \in \mathbb{R}^{k-1}$ satisfying
$0\le y_1 \le \cdots \le y_{k-1} \le J$ and $j_i-1 < y_i \le j_i$ for each $i$
(there are $b_j$ numbers $y_i$ in each interval $(j-1,j]$).
Making the change of variables $\xi_i=y_i/J$ and summing over all possible
vectors $(j_1,\ldots,j_{k-1})\in \mathcal{J}$, 
we find that the inner sum in \eqref{Nku1} is
\begin{align*}
&\le (2J)^{k-1} \Vol \left\{ 0\le \xi_1 \le \cdots \le \xi_{k-1} \le 1 : \xi_i
\ge \frac{\a i - \b - K'-A-2}{2J} \; (1\le i\le k-1) \right\} \\
&\le \frac{(\log_2 x+2K')^{k-1}}{(k-1)!} Q_{k-1} \( \frac{\b+K'+A+2}{\a}, 
\frac{2J}{\a} \) \\
&\ll_A \frac{(\log_2 x)^{k-1}}{(k-1)!} \; \frac{(u+1)w}{k},
\end{align*}
where we have used \eqref{Q}.  By \eqref{Nku1}, summing on 
$\ell$ and using \eqref{pik} completes the proof.

%
\section{Proof of Theorem \ref{Nk} lower bound}
%

First, we assume $k\ge 2$, since if $k=1$ then $N_1(x;\a,\b) = \pi_1(x) +
O(\log x)$
trivially as $A+\b \ge \a$ (powers of primes $\le e^{\a-\b}$ are not counted
in $N_1(x;\a,\b)$).  
Also, we may assume that $\a \ge 1/2A$.  If $\a < 1/2A$, then
$N_k(x;\a,\b) \ge N_k(x;1/2A,0)$ and we prove below that
$N_k(x;1/2A,0)\gg \pi_k(x)$ (here $u=0$, $v\ge 2k$ and $w\ge k$).

Let $T$ be a sufficiently large
constant, depending on $\eps$ and $A$, and put 
$$
C = e^{3T + 2K +10}.
$$
We first prove the theorem in the case that 
\be\label{Nkass}
e^{\a(w-1)} - e^{\a(w-2)} \ge C.
\ee
Notice that 
\be\label{ajbeta}
\a j - \b = \log_2 x - \a (w+k-1-j).
\ee
In particular,
$$
\a k - \b = \log_2 x - \a (w-1) \le \log_2 x - \log C.
$$
Let $J=\fl{\log_2 x - K - \log T - 2}$.
Recall the definition of the numbers $\lambda_j$ and sets $G_j$ from
section \ref{partitions}.
Consider squarefree $n$ satisfying \eqref{pjlower},
with $p_{k-1}\le \lambda_J$ and for which
$$
p_1 \cdots p_{k-1} \le x^{1/2}.
$$
Also take $p_k$ so that $x/2 < n \le x$.  Given $p_1,\ldots,p_{k-1}$, the
number of possible $p_k$ is $\gg x/(p_1\cdots p_{k-1} \log x)$.
Put $b_1=\cdots=b_{T-1}=0$ and for $T \le j\le J$, suppose
$b_j \le \min (T(j-T-1),T(J-j+1))$.  Suppose there are
exactly $b_j$ primes $p_i$ in the set $G_j$ for $1\le j\le J$.  By the
definition of $J$,
$$
\sum_{i=1}^{k-1} \log p_i \le T e^{J+K} \sum_{r=1}^{k-1} r e^{1-r}
< 3 T e^{J+K} \le \frac12 \log x, 
$$
as required.  Define the numbers $j_i$ by $p_i \in G_{j_i}$.
The inequalities \eqref{pjlower} will be satisfied
if 
\be\label{Nklj}
j_i \ge \a i - \b + K \qquad (1\le i\le k-1).
\ee
This is possible since by \eqref{Nkass},
$$
\a(k-1)-\b = \log_2 x - \a w \le  \log_2 x - 2K - 3T - 10 < J-T-1.
$$
With $(j_1,\ldots,j_{k-1})$ fixed (so that $b_1,\ldots,b_J$ are fixed),
the sum of $1/p_1\cdots p_{k-1}$ is
\begin{align*}
&= \prod_{j=T}^J \frac{1}{b_j!} \biggl(
  \sum_{p_1\in G_j} \frac{1}{p_1} \sum_{\substack{p_2\in G_j \\ p_2 \ne p_1}}
  \frac{1}{p_2} \cdots \sum_{\substack{p_{b_j}\in G_j \\ p_{b_j} \not\in
  \{ p_1, \ldots, p_{b_j-1} \} }} \frac{1}{p_{b_j}} \biggr) \\
&\ge\prod_{j=T}^J \frac{1}{b_j!}
   \biggl( \sum_{p\in G_j} \frac{1}{p} - \frac{b_j-1}{\lam_{j-1}}
   \biggr)^{b_j} \\
&\ge \prod_{j=T}^J \frac{1}{b_j!}
   \biggl( 1 - \frac{b_j}{\lam_{j-1}} \biggr)^{b_j} \\
&\ge \prod_{j=T}^J \frac{1}{b_j!}  \biggl( 1 - \frac{T(j-T+1)}{\exp\exp(j-1-K)}
   \biggr)^{T(j-T+1)} \\
&\ge \frac{1/2}{b_T! \cdots b_J!}
\end{align*}
if $T$ is large enough.  The right side is $1/2$ of the volume of the region
of $(y_1,\cdots,y_{k-1}) \in \mathbb{R}^{k-1}$ satisfying
$0\le y_1 \le \cdots \le y_{k-1} \le J-T+1$ and $j_i-T \le y_i \le j_i+1-T$
for each $i$.  Set $H=J-T+1$.  Assume that 
\be\label{jmT}
j_{mT+1} \ge T+m, \quad j_{k-1-mT} \le J-m \qquad (\text{integers } m\ge 1),
\ee
so that $b_j \le \min (T(j-T+1),T(J-j+1))$ for each $j$.
Making the substitution $\xi_i = y_i/H$ and summing
over all tuples $(j_1,\cdots,j_{k-1})$ yields
\be\label{Nkl2}
N_k(x;\a,\b) \gg \frac{x H^k}{\log x} \text{Vol} (R) \gg_A \frac{x(\log_2
  x)^k}{\log x} \text{Vol} (R),
\ee
where, by \eqref{Nklj} and \eqref{jmT}, $R$ is the set of $\bxi$ satisfying
(i) $0\le \xi_1 \le \cdots \le \xi_{k-1} \le 1$,
$\xi_i \ge (\a i - \b + K - T)/H$ for each $i$, 
(ii) $\xi_{mT+1} \ge m/H$ and 
$\xi_{k-1-mT} \le 1 - m/H$ for each positive integer
$m$.  

It remains to estimate from below the volume of $R$.
Let $S$ be the set of $\bxi$ satisfying (i), so that
$$
\text{Vol} (S) = \frac{Q_{k-1}(\mu,\nu)}{(k-1)!}, \qquad \mu=\frac{\b+T-K}{\a},
\; \nu = \frac{H}{\a}.
$$
If $T\ge K+A$, then $\mu \asymp_A u+1$.  By the definition of $C$ and $J$, 
if $T$ is large enough then
$$
\mu+\nu-(k-1) = \frac{J-K+1+\b}{\a}-(k-1) \ge w - \frac{\log T+2K+2}{\a}
\ge \frac{w}{1+\eps} \ge 1.
$$
Hence, by \eqref{Q}, 
\be\label{volS}
\text{Vol} (S) \gg \frac{f}{(k-1)!}, \qquad f = \min(1,(u+1)w/k). 
\ee
The implied constant in \eqref{volS} does not depend on $T$, 
but the inequality does require that $T$ be sufficiently large.

For a positive integer $m$, let 
\begin{align*}
V_1(m) &= \text{Vol} \{ \bxi \in S : \xi_{mT+1} < m/H \}, \\
V_2(m) &= \text{Vol} \{ \bxi \in S : \xi_{k-1-mT} > 1- m/H\}.
\end{align*}
We have by \eqref{Q},
\begin{align*}
V_1(m) &\le \frac{(m/H)^{mT+1}}{(mT+1)!} \text{Vol} \{ 0\le \xi_{mT+2} \le
  \cdots \le \xi_{k-1} \le 1 : \xi_i \ge \tfrac{i-\mu}{\nu} \; 
  (mT+2\le i\le k-1) \} \\
&= \frac{(m/H)^{mT+1}}{(mT+1)!} \;
  \frac{Q_{k-2-mT}(\mu-(mT+1),\nu)}{(k-2-mT)!} \\
&\ll \frac{(m/H)^{mT+1}}{(mT+1)!} \frac{\mu(\mu+\nu-(k-1))}{(k-mT)(k-2-mT)!}\\
&\ll \frac{ f k (m/H)^{mT+1}}{(k-mT)(mT+1)!(k-2-mT)!} \\
&\le \frac{f}{(k-1)!} \; \frac{(km/H)^{mT+1}}{(mT+1)!} \; \frac{k}{k-mT}.
\end{align*}
Since $k/H \ll_A 1$ and $r! \ge (r/e)^r$, it follows from \eqref{volS} 
that for large enough $T$,
$$
\sum_m V_1(m) \le \frac14 \text{Vol} (S).
$$
Similarly,
$$
V_2(m) \le \frac{Q_{k-2-mT}(\mu,\nu)}{(k-2-mT)!} \frac{(m/H)^{mT+1}}{(mT+1)!}.
$$
By \eqref{Q},
$$
Q_{k-2-mT}(\mu,\nu) \ll \min \( 1, \frac{\mu(\mu+\nu-(k-1)+mT+1)}{k-mT} \) 
\ll \frac{mT k f}{k-mT}.
$$
Hence, if $T$ is large enough then
$$
\sum_m V_2(m) \le \frac14 \text{Vol} (S).
$$
We therefore have, for $T$ large enough,
$$
\text{Vol} (R) \ge \text{Vol} (S) - \sum_{m\ge 1} (V_1(m)+V_2(m)) \gg_A
\frac{f}{(k-1)!}.
$$
Together with \eqref{Nkl2} and \eqref{pik}, this completes the proof
under the assumption \eqref{Nkass}.

It remains to consider the case 
$$
1+\eps \le e^{\a(w-1)} - e^{\a(w-2)} \le C.
$$
Since $w\ge 1+\eps$ and $\a \ge 1/2A$, 
we find that $\a \ll_{\eps,A} 1$ and $w\ll_{\eps,A} 1$.  Hence, if
$x$ is large enough,
$$
k = u+v-w+1 \ge v-w \ge \frac{\log_2 x}{4A}.
$$
Let $B$ be a large integer depending on $\eps$.  Suppose that
\be\label{largej}
\a j - \b \le \log_2 p_j \le \a j-\b + \log(1+\eps/2) \qquad (k-B \le j\le k-1)
\ee
Then, by \eqref{ajbeta},
\begin{align*}
\sum_{j=k-B}^{k-1} \log p_j &\le (1+\eps/2)\( e^{-\a w} + e^{-\a(w+1)} + \cdots
+ e^{-\a(w+B-1)} \) \log x 
\\ &< (1+\eps/2) \( \frac{1}{e^{\a(w-1)}-e^{\a(w-2)}} -e^{-\a(w-1)} \) \log x.
\end{align*}
Assume also that 
\be\label{smallj}
\sum_{j=1}^{k-B-1} \log p_j \le \frac{\eps/2} {e^{\a(w-1)}-e^{\a(w-2)}} \log x.
\ee
If in addition $\a k - \b \le \log_2 p_k \le \a k - \b + \log(1+\eps/2)$,
then by \eqref{exp1},
$$
\log n = \sum_{j=1}^k \log p_j \le \frac{\eps/2 + 1 + \eps/2} 
{e^{\a(w-1)}-e^{\a(w-2)}}\log x \le \log x,
$$
as required.  Thus, given $p_1,\ldots,p_{k-1}$ satisfying
\eqref{largej} and \eqref{smallj}, the number of $p_k$ is
$\gg x/(p_1 \cdots p_{k-1}\log x)$.  If $B$ is large enough,
there is great flexibility in choosing $p_1,\ldots,p_{k-B-1}$, 
since by \eqref{ajbeta},
$$
\sum_{j=1}^{k-B-1} e^{\a j - \b} \le \frac{e^{-\a(B+1)}}
{e^{\a(w-1)}-e^{\a(w-2)}} \log x,
$$
which is small compared with the right side of \eqref{smallj}.
By the same argument used to give  a lower bound for the sum of
$1/(p_1\cdots p_{k-1})$ under the assumption \eqref{Nkass}, we obtain
$$
\sum_{p_1,\ldots,p_{k-B-1}} \frac{1}{p_1 \cdots p_{k-B-1}} \gg_{A,\eps}
\frac{f (\log_2 x)^{k-B-1}}{(k-B-1)!}.
$$
Also, since $k\gg_A \log_2 x$, we have
$$
\sum_{p_{k-B},\ldots,p_{k-1}} \frac{1}{p_{k-B}\cdots p_{k-1}} \gg_{\eps,B}
1 \gg_{\eps,A} (\log_2 x)^{B} \frac{(k-B-1)!}{(k-1)!}.
$$
The proof is again completed by applying \eqref{pik}.

%
%

\bibliographystyle{plain}
\bibliography{wnt}

\end{document}